\documentclass{amsart}
\usepackage{amsmath}
\usepackage{amstext}
\usepackage{amssymb}
\usepackage{amsthm}
\swapnumbers
\theoremstyle{plain}
\newtheorem{thm}{Theorem}[section]
\newtheorem{lem}[thm]{Lemma}

\newtheorem{prop}[thm]{Proposition}

\theoremstyle{definition}

\newtheorem{defrmks}[thm]{Definition and Remarks}

\newtheorem{ntn}[thm]{Notation}

\theoremstyle{remark}

\newtheorem{rmk}[thm]{Remark}

 \DeclareMathOperator{\Hom}{Hom}
 
 \DeclareMathOperator{\Ker}{Ker}

\DeclareMathOperator{\id}{inj\,dim}
\DeclareMathOperator{\rk}{rank}

\def\Z{\mathbb Z}
\def\N{\mathbb N}

\def\fa{{\mathfrak{a}}}

\def\fm{{\mathfrak{m}}}

\def\fn{{\mathfrak{n}}}

\def\nn{\relax\ifmmode{\mathbb N_{0}}\else$\mathbb N_{0}$\fi}
\def\lra{\longrightarrow}

\begin{document}

\title[ON THE HARTSHORNE--SPEISER--LYUBEZNIK THEOREM]{ON THE
HARTSHORNE--SPEISER--LYUBEZNIK THEOREM
ABOUT ARTINIAN MODULES WITH A FROBENIUS ACTION}
\author{RODNEY Y. SHARP}
\address{Department of Pure Mathematics,
University of Sheffield, Hicks Building, Sheffield S3 7RH, United Kingdom\\
{\it Fax number}: 0044-114-222-3769}
\email{R.Y.Sharp@sheffield.ac.uk}

\thanks{The author was partially supported by the
Engineering and Physical Sciences Research Council of the United
Kingdom (grant number EP/C538803/1).}


\date{\today}


\begin{abstract}
Let $R$ be a commutative Noetherian local ring of prime
characteristic. The purpose of this paper is to provide a short
proof of G. Lyubeznik's extension of
a result of R. Hartshorne and R. Speiser about a module over the skew polynomial
ring $R[x,f]$ (associated to $R$ and the Frobenius homomorphism
$f$, in the indeterminate $x$) that is both $x$-torsion and
Artinian over $R$.
\end{abstract}

\maketitle

\setcounter{section}{-1}
\section{\sc Introduction}
\label{in}

In the theory of tight closure of ideals in a $d$-dimensional
commutative (Noetherian)
local ring $(R,\fm)$ of prime characteristic $p$, study of properties of the
`top' local cohomology module $H^d_{\fm}(R)$ related to the Frobenius
homomorphism $f: R \lra R$ has been a very effective tool: see, for example,
K. E. Smith \cite{S94,S95}. Some of the properties of
$H^d_{\fm}(R)$ related to $f$ can be neatly described in terms of a natural
structure which $H^d_{\fm}(R)$ possesses as a left module over the
skew polynomial ring $R[x,f]$; also, it is well known that $H^d_{\fm}(R)$
is Artinian as an $R$-module. One can take the view that $H^d_{\fm}(R)$
is an important example of a left $R[x,f]$-module that is Artinian as an $R$-module.

In 1977, R. Hartshorne and R. Speiser \cite[Proposition
1.11]{HarSpe77} proved, in the case where the local ring $R$ of
characteristic $p$ contains its residue field which is perfect,
that, given a left $R[x,f]$-module $H$ that is Artinian as an $R$-module,
there exists a non-negative integer $e$ with the following property: whenever
$h \in H$ is such that $x^jh = 0$ for some positive integer $j$, then $x^eh = 0$.

Twenty years later, G.
Lyubeznik \cite[Proposition 4.4]{Lyube97} proved this result
without restriction on the local ring $R$ of characteristic $p$,
that is, he was able to drop the hypotheses about the residue
field of $R$. Lyubeznik's proof is an application of his
substantial theory of $F$-modules.

There is some evidence that the Hartshorne--Speiser--Lyubeznik
Theorem can be exploited to good effect in tight closure theory.
For example, it has recently been used in \cite{tctecpi} to prove that, if $c$
is a test element for a reduced excellent equidimensional local ring $(R,\fm)$
of characteristic $p$, then there exists
a power of $p$ that is a test exponent for $c, \fa$ (see
\cite[Definition 2.2]{HocHun00}) for every parameter ideal $\fa$ of $R$
simultaneously.

It therefore seems desirable to
have a short proof of the Hartshorne--Speiser--Lyubeznik
Theorem that does not rely on the theory
of $F$-modules. This paper provides one such that actually
follows the general line of the
Hartshorne--Speiser proof.

\section{\sc Left modules over the skew polynomial ring $R[x,f]$}
\label{sp}

\begin{ntn}
\label{nt.1} Throughout the paper, $A$ will denote a general commutative
Noetherian ring and $R$ will denote a commutative
Noetherian ring of prime characteristic $p$. In cases where
such a ring is assumed to be local, the notation $(A,\fm)$ or
$(R,\fm)$ will indicate that $\fm$ is the maximal ideal.

We shall always denote by $f:R\lra R$ the Frobenius homomorphism,
for which $f(r) = r^p$ for all $r \in R$.  We use $\N$ and $\nn$
to denote the sets of positive integers and non-negative integers,
respectively.  We shall work with the
 skew polynomial ring $R[x,f]$ associated to $R$ and $f$
in the indeterminate $x$ over $R$. Recall that $R[x,f]$ is, as a
left $R$-module, freely generated by $(x^i)_{i \in \nn}$,
 and so consists
 of all polynomials $\sum_{i = 0}^n r_i x^i$, where  $n \in \nn$
 and  $r_0,\ldots,r_n \in R$; however, its multiplication is subject to the
 rule
 $$
  xr = f(r)x = r^px \quad \mbox{~for all~} r \in R\/.
 $$
\end{ntn}

\begin{defrmks}
\label{nt.2}
We say that the left $R[x,f]$-module $H$ is {\em $x$-torsion-free\/} if $xh = 0$, for
$h \in H$, only when $h = 0$.
The set $\Gamma_x(H) := \left\{ h \in H : x^jh = 0
\mbox{~for some~} j \in \N \right\}$ is an $R[x,f]$-submodule of
$H$, called the\/ {\em $x$-torsion submodule} of $H$. In general, the
$R[x,f]$-module $H/\Gamma_x(H)$ is $x$-torsion-free.
\end{defrmks}

\section{\sc The Hartshorne--Speiser Theorem}
\label{hs}

As explained in the Introduction, this paper is concerned
with the following result of R. Hartshorne and
R. Speiser.

\begin{thm}[Hartshorne--Speiser {\cite[Proposition
1.11]{HarSpe77}}] \label{hs.1} Suppose that $R$ is local and
contains its residue field which is perfect. Let $H$ be a left
$R[x,f]$-module which is Artinian as an $R$-module. Then there
exists $e \in \nn$ such that $x^e\Gamma_x(H) = 0$.
\end{thm}

G. Lyubeznik \cite[Proposition 4.4]{Lyube97} proved this result
without restriction on the local ring $R$ of characteristic $p$,
that is, he was able to drop the hypotheses about the residue
field of $R$; his proof is an application of his theory of
$F$-modules. The main purpose of this section is to show how one
can modify the argument of Hartshorne and Speiser to obtain a
short and direct proof of the result in the generality achieved by
Lyubeznik. To
achieve this aim, we shall establish a generalization of
Proposition 1.9 of Hartshorne--Speiser \cite{HarSpe77}.

Our first preparatory result concerns an Artinian module of finite
injective dimension over a general local ring $(A,\fm)$. Let $E$
denote $E_A(A/\fm)$, the injective envelope of the simple
$A$-module $A/\fm$. Recall that an $A$-module is Artinian if and
only if it is isomorphic to a submodule of $E^t$, the direct sum
of $t$ copies of $E$, for some $t \in \N$. It follows that, if $G$
is an Artinian $A$-module, then, for each $i \in \nn$, the $i$-th
term $E^i_A(G)$ in the minimal injective resolution of $G$ is
isomorphic to a direct sum of finitely many copies of $E$. When
$J$ is an Artinian injective $A$-module, we shall use the Bass
number $\mu^0(\fm,J)$ to denote the number of copies of $E$ that
occur in a decomposition of $J$ as a direct sum of indecomposable
injective $A$-modules.

\begin{prop}
\label{hs.1a} Let $G$ be an Artinian module over the local ring
$(A,\fm)$ such that $\id _AG < \infty$.

\begin{enumerate}

\item Let
$$
I^{\scriptscriptstyle \bullet} :
0 \longrightarrow
I^0 \stackrel{d^0}{\longrightarrow}
I^1 \longrightarrow \cdots \longrightarrow
I^i \stackrel{d^i}{\longrightarrow}
I^{i+1} \longrightarrow \cdots
$$
be a finite injective resolution of $G$ in which each term is
isomorphic to a direct sum of copies of $E := E_A(A/\fm)$. (It
should be noted that the minimal injective resolution of $G$ has
this property.) Then the integer
$\sum_{i=0}^{\infty}(-1)^i\mu^0(\fm,I^i)$ is independent of the
choice of finite injective resolution $I^{\scriptscriptstyle
\bullet}$ of $G$ having the stated property. We call this integer
the\/ {\em Euler number} of $G$, and denote it by $\chi(G)$ (or
$\chi_A(G)$ when it is desirable to emphasize the local ring $A$).

\item The Euler number $\chi(G)$ of $G$ is non-negative.

\item Let $G', \overline{G}$ be further Artinian $A$-modules of
finite injective dimension and
suppose that there is an exact sequence $0 \lra G' \lra G \lra
\overline{G} \lra 0$ in the category of $A$-modules and $A$-homomorphisms. Then
$$
\chi (G) = \chi (G') + \chi (\overline{G}).
$$

\item When $A$ is complete, the following three conditions are equivalent:
\begin{enumerate}
\item $(0:_AG) \neq 0$;
\item $\chi(G) = 0$;
\item $(0:_AG)$ contains a non-zerodivisor of $A$.
\end{enumerate}

\end{enumerate}
\end{prop}

\begin{proof}
(i),(ii),(iii) There is an $A$-homomorphism $\alpha : G \longrightarrow
I^0$ such that the sequence
\[
0
\longrightarrow
G \stackrel{\alpha}{\longrightarrow}
I^0 \stackrel{d^0}{\longrightarrow}
I^1 \longrightarrow \cdots \longrightarrow
I^i \stackrel{d^i}{\longrightarrow}
I^{i+1} \longrightarrow \cdots
\]
is exact. Note that $E$, $G$, $G'$, $\overline{G}$ and all the $I^j~(j \in \nn)$ have
natural structures as modules over the completion
$(\widehat{A},\widehat{\fm})$ of $A$, and that, when they are
given these, there is an $\widehat{A}$-isomorphism $E \cong
E_{\widehat{A}}(\widehat{A}/\widehat{\fm})$ and
the above-displayed exact sequence provides an
injective resolution of $G$ as an $\widehat{A}$-module. Furthermore,
$0 \lra G' \lra G \lra
\overline{G} \lra 0$ is an exact sequence
in the category of $\widehat{A}$-modules and $\widehat{A}$-homomorphisms.
It thus follows that
it is sufficient to prove parts (i), (ii) and (iii) under the additional
assumption that $A$ is complete.

Let $D$ be the functor $\Hom_{A}(\: {\scriptscriptstyle \bullet}
\:,E)$ on the category of $A$-modules. We use Matlis duality.
Since $D(E) \cong A$, application of the functor $D$ to
$I^{\scriptscriptstyle \bullet}$ yields an exact sequence
$$
\cdots \longrightarrow D(I^{i+1}) \longrightarrow D(I^i) \lra \cdots \lra
D(I^0) \lra D(G) \lra 0,
$$
and this provides a finite free resolution of the finitely
generated $A$-module $D(G)$. Moreover, for each $i \in \nn$, the
free $A$-module $D(I^i)$ is finitely generated of rank
$\mu^0(\fm,I^i)$. Thus
$$
\sum_{i=0}^{\infty}(-1)^i\mu^0(\fm,I^i) =
\sum_{i=0}^{\infty}(-1)^i\rk D(I^i),
$$
which is just the Euler number $\chi(D(G))$, and so is independent of the
choice of finite injective resolution $I^{\scriptscriptstyle \bullet}$ of $G$ of
the type under consideration: see \cite[pp.\ 159]{HM}, for
example. Likewise, the claim in part (ii) now follows from
the corresponding statement (see \cite[Theorem 19.7]{HM}, for example)
about modules with finite free resolutions, and the claim in part (iii)
follows from the well-known fact that
$\chi$ is additive on short exact sequences of modules with finite free resolutions.

(iv) Since the annihilators of $G$ and $D(G)$ are equal, the
equivalence of (a), (b) and (c) is now immediate from a theorem
of M. Auslander and D. A. Buchsbaum \cite{AusBuc58}: see
\cite[Theorem 19.8]{HM}, for example.
\end{proof}

\begin{rmk}
\label{hs.2} It is a consequence of Proposition \ref{hs.1a} that,
with the notation of that result, $\chi_A(G) =
\sum_{i=0}^{\infty}(-1)^i\mu^i(\fm,G) = \sum_{i=0}^{\id
G}(-1)^i\mu^i(\fm,G)$, because, for each $i \in \nn$, the $i$-th
term in the minimal injective resolution of $G$ is isomorphic to
the direct sum of $\mu^i(\fm,G)$ copies of $E$, and $\mu^j(\fm,G)
= 0$ for all $j > \id_A G$.
\end{rmk}

We can now establish the promised generalization of Proposition
1.9 of Hartshorne--Speiser \cite{HarSpe77}.

\begin{prop}
\label{hs.3} {\rm (Compare Hartshorne--Speiser \cite[Proposition
1.9]{HarSpe77}.)} Assume that $(R,\fm)$ is a complete regular
local ring, and that $H$ is a left $R[x,f]$-module which is
Artinian as an $R$-module and such that $RxH = H$. Let $K = \{ h
\in H : xh = 0\}$, an $R[x,f]$-submodule of $H$. Then $(0:_RK)
\neq 0$.
\end{prop}

\begin{proof} Here, we shall use $R'$ to denote $R$ considered
as an $R$-module by means of
$f$ (at points where care is needed). Also $F$ will denote the
{\em Frobenius functor\/} $R'\otimes _R(\: {\scriptscriptstyle
\bullet} \:)$ from the category of all $R$-modules and
$R$-homomorphisms to the category of all $R'$-modules and
$R'$-homomorphisms.

Since $axrh = ar^pxh$ for $a \in R'$, $r \in R$ and $h \in H$,
there is an $R$-homomorphism $\phi : F(H) \lra H$ for which
$\phi(a\otimes h) = axh$ for all $h \in H$ and $a \in R'$. Note
that $\phi$ is surjective, because $RxH = H$. Note also that, if
$h \in K$, then the element $1 \otimes h$ of $F(H)$ lies in $\Ker
\phi$. Since $R$ is regular, $f : R \lra R$ is flat (by E. Kunz
\cite{Kunz69}), and therefore faithfully flat. The
$\Z$-homomorphism $\gamma : K \lra \Ker \phi$ for which $\gamma
(h) = 1 \otimes h$ for all $h \in K$ is therefore injective. It is
therefore enough for us to show that $(0:_{R}\Ker \phi) \neq 0$,
for if $0 \neq a \in R$ annihilates $\Ker \phi$, then, for each $h
\in K$, we have $\gamma(ah) = 1 \otimes ah = a^p \otimes h = a^p
(1 \otimes h) = 0$.

There is a short exact sequence
$$
0 \lra \Ker \phi \hookrightarrow F(H) = R'\otimes _RH
\stackrel{\phi}{\lra} H \lra 0$$ of $R$-modules and
$R$-homomorphisms. By Huneke--Sharp \cite[Proposition 1.5]{58},
for each injective $R$-module $I$, we have $F(I) \cong I$. Observe
that every $R$-module has finite injective dimension, because $R$
has finite global dimension. If one applies the exact functor $F$
to the minimal injective resolution for $H$, one can deduce, with
the aid of Proposition \ref{hs.1a}, that $F(H)$ is isomorphic to a
submodule of the direct sum of finitely many copies of $F(E) \cong
E$ and so is Artinian, and that $\chi (F(H)) = \chi(H)$. Hence
$\Ker \phi$ is an Artinian $R$-module, and it follows from
Proposition \ref{hs.1a}(iii) that $\chi (\Ker \phi) = \chi (F(H)) -
\chi(H) = 0$. Hence $(0:_{R}\Ker \phi) \neq 0$ by
Proposition \ref{hs.1a}(iv).
\end{proof}

We shall need the following lemma of Hartshorne and Speiser.

\begin{lem} [Hartshorne--Speiser {\cite[Lemma 1.10]{HarSpe77}}]
\label{hs.3p} Let $H$ be a left $R[x,f]$-module, and set $$K := \{
h \in H : xh = 0\},$$ an $R[x,f]$-submodule of $H$. Suppose that
$a \in R$ is such that $aK = 0$. Then $a^2\Gamma_x(H) = 0$.
\end{lem}

The short proof, presented in the next theorem, of
Lyubeznik's extension of the Hartshorne--Speiser Theorem follows
the general line of argument of Hartshorne and Speiser.

\begin{thm} [G. Lyubeznik {\cite[Proposition 4.4]{Lyube97}}]
\label{hs.4} {\rm (Compare
Hartshorne--Speiser \cite[Proposition 1.11]{HarSpe77}.)}
Suppose that $(R,\fm)$ is local, and let $H$ be a
left $R[x,f]$-module which is Artinian as an $R$-module.
Then there
exists $e \in \nn$ such that $x^e\Gamma_x(H) = 0$.
\end{thm}

\begin{proof}
Recall the natural $\widehat{R}$-module structure on the Artinian
$R$-module $H$: given $h \in H$, there exists $t \in \N$ such that
$\fm^th = 0$; for an $\widehat{r} \in \widehat{R}$, choose any $r
\in R$ such that $\widehat{r} - r \in \fm^t\widehat{R}$; then
$\widehat{r}h = rh$. It is easy to see from this that
$x\widehat{r}h = \widehat{r}^pxh$ for all $h \in H$ and
$\widehat{r} \in \widehat{R}$, and we can then use \cite[Lemma
1.3]{KS} to see that one can assume that $R$ is complete.

Argue by induction on $n := \dim R$; note that, when $n = 0$, the
Artinian $R$-module $H$ has finite length and then the claim
follows easily. Suppose that $n > 0$ and assume
inductively that the result has been proved when the underlying
complete local ring $R$ has dimension smaller than $n$.

Let $k$ be a coefficient field for $R$, and let $r_1, \ldots, r_n$
be a system of parameters for $R$. Then $R$ is module-finite over
the complete regular local ring $k[[r_1, \ldots, r_n]]$, which we
denote by $(S,\fn)$. Since $\fm \cap S = \fn$, it is clear that
each element of $H$ is annihilated by some power of $\fn$. Since
$\fn R$ is $\fm$-primary, it contains $\fm^t$ for some $t \in \N$,
and so $(0:_H \fn) = (0:_H \fn R) \subseteq (0:_H \fm^t)$, which
is finitely generated over $R$ and therefore over $S$. Thus $H$ is
Artinian as $S$-module.

We can replace $H$ by its $R[x,f]$-
and $S[x,f]$-submodule $\Gamma_x(H)$; thus we can
assume that $H$ is $x$-torsion. Of course, we can assume that $H
\neq 0$.

The descending chain of $S[x,f]$-submodules
$$
H \supseteq SxH \supseteq Sx^2H \supseteq \cdots \supseteq Sx^iH
\supseteq Sx^{i+1}H \supseteq \cdots
$$
of $H$ must eventually stabilize: let $t \in \nn$ be such that
$Sx^tH = Sx^{t+j}H$ for all $j \in \N$. Observe that $Sx^tH =
Sx(Sx^tH)$, and that it is enough to prove the claim
for $Sx^tH$ rather than $H$.

It thus follows that, in order to complete the inductive step, we can (replace $R$
by $S$ and) assume that $R$ is regular (and complete), that $H$ is
$x$-torsion and that $H = RxH$.

Let $K = \{ h \in H : xh = 0\}$. By Proposition \ref{hs.3}, there
exists $0 \neq a \in R$ such that $aK = 0$.
Therefore $a^2\Gamma_x(H) = a^2 H = 0$, by Lemma \ref{hs.3p}.

Thus $H$ has a natural structure as a module over the complete
local ring $R/Ra^2$, which has dimension $n-1$. Use $\overline{r}$
to denote the natural image in $R/Ra^2$ of an element $r \in R$.
Then $x\overline{r}h = \overline{r}^pxh$ for all $r \in R$ and $h
\in H$; it thus follows from \cite[Lemma 1.3]{KS} that $H$
inherits a structure as left $(R/Ra^2)[x,f]$-module, compatible
with its $R[x,f]$-module structure; note that $H$ is still
Artinian over $R/Ra^2$ and $x$-torsion, and satisfies $H =
(R/Ra^2)xH$. Application of the inductive hypothesis therefore
completes the proof.
\end{proof}

I am grateful to Craig Huneke for pointing out that
the argument used in the above proof can be modified to prove
the next result.

\begin{thm}
\label{hs.5}
Suppose that $(R,\fm)$ is local, and let $H$ be an $x$-torsion
left $R[x,f]$-module which is Artinian as an $R$-module.
If $H = RxH$, then $H$ has finite length as an $R$-module.
\end{thm}

\begin{proof} Only a sketch is presented here, as the strategy used is very
similar to that used in the above proof of Theorem \ref{hs.4}.

One can assume that $R$ is complete; then argue by induction on $n
:= \dim R$, the result being easy when $n = 0$. For the inductive
step, in the situation where $n >0$, again let $k$ be a
coefficient field for $R$, let $r_1, \ldots, r_n$ be a system of
parameters for $R$, and let $S$ be the complete regular local ring
$k[[r_1, \ldots, r_n]]$. We can again use the fact that $R$ is
module-finite over $S$ to see that $H$ is Artinian as $S$-module,
so that there exists $t \in \nn$ such that $Sx^tH = Sx^{t+j}H$ for
all $j \in \N$.

Thus the left $S[x,f]$-module $Sx^tH$ is $x$-torsion, Artinian as
an $S$-module, and such that $Sx^tH = Sx(Sx^tH)$. Let $e_1 = 1,
e_2, \ldots, e_h$ generate $R$ as an $S$-module. If we could prove
that $Sx^tH$ is of finite length as an $S$-module, then, since $H
= RxH$, it would follow that
$$
H = Rx^tH = \sum_{i=1}^h Se_ix^tH = \sum_{i=1}^h e_iSx^tH
$$
is finitely generated as an $S$-module, and therefore as an
$R$-module. It follows that, in order to complete the inductive step, we can
(replace $R$
by $S$ and) assume that $R$ is regular (and complete). The proof can now be
completed by an argument almost identical to that in the last two paragraphs of
the above proof of Theorem \ref{hs.4}.
\end{proof}

\begin{rmk}
\label{hs.6} Note that an extension of the Hartshorne--Speiser--Lyubeznik Theorem
to non-local situations is proved in \cite[Corollary 1.8]{tctecpi}: there
it is proved that, if $R$ is merely a commutative Noetherian ring (of
characteristic $p$), and
$H$ is a left $R[x,f]$-module which is Artinian as an $R$-module, then
there exists $e \in \nn$ such that $x^e\Gamma_x(H) = 0$.
\end{rmk}

\bibliographystyle{amsplain}

\end{document}